\documentclass{article}
\usepackage{amssymb}
\title{{\bf On a  Class Number Formula for Real Quadratic
             Number Fields}\\ }
\author{David M.~Bradley, Ali E.~\"{O}zl\"{u}k, C.~Snyder} 
\date{July 19, 2001}

\begin{document}
\maketitle
\newtheorem{thm}{Theorem}
\newtheorem{Def}{Definition}
\newtheorem{lem}{Lemma}
\newtheorem{cor}{Corollary}
\newtheorem{prop}{Proposition}
\newcommand{\Q}{{\mathbb Q}}
\newcommand{\N}{{\mathbb N}}
\newcommand{\Z}{{\mathbb Z}}
\newcommand{\C}{{\mathbb C}}
\newcommand{\eps}{\varepsilon}
\newcommand{\eop}{\hfill $\square$}
\newcommand{\R}{{\rm Re}}
\newcommand{\I}{{\rm Im}}
\newcommand{\tfrac}[2]{\textstyle\frac{#1}{#2}\displaystyle}

\begin{abstract}
For an even Dirichlet character $\psi$, we obtain a formula for
$L(1,\psi)$ in terms of a sum of Dirichlet $L$-series evaluated at
$s=2$ and $s=3$ and a rapidly convergent numerical series
involving the central binomial coefficients. We then derive a
class number formula for real quadratic number fields by taking
$L(s,\psi)$ to be the quadratic $L$-series associated with these
fields.
\end{abstract}

\section{Introduction}
In~\cite{B}, acceleration formul{\ae} are derived for Catalan's
constant $L(2,\chi_4)$; (here $\chi_4$ is the non-principal
Dirichlet character of modulus $4$). In some of these formul{\ae}
$L(2,\chi_4)$ is given as the sum of two terms: one involving a
rapidly convergent series and the other involving the natural
logarithm of a unit in the ring of integers of a finite abelian
field extension of the rational number field $\Q$.  The existence
of the logarithmic terms suggested to the authors that these terms
should somehow be related to the values of Dirichlet $L$-series at
the argument $s=1$. This leads to the general question of whether
or not there exist relations between the value of $L$-series at
$s=1$ and values of $L$-series at integer arguments larger than
$1$.

The purpose of this note is to exhibit such a relation between
values of $L$-series. For an even Dirichlet character $\psi$, we
obtain a formula for $L(1,\psi)$ in terms of a sum of Dirichlet
series evaluated at $s=2$ and $s=3$ and a convergent numerical
series involving powers of twice special values of the sine
function divided by ${2n\choose n}n^3$. See Theorem~\ref{T} below
for a precise statement.  (It is perhaps interesting to notice
that not much is known about number theoretic properties of the
values of the $L$-series on the right-hand side of the formula
given in this theorem.)  We then deduce a class number formula for
real quadratic number fields by letting $\psi$  be the quadratic
character associated with a real quadratic number field; see
Corollary~\ref{C}.  This class number formula seems new to us and
is perhaps an interesting curiosity.

To derive our results, we employ a formula of Zucker~\cite{Z} that
expresses
\begin{equation}
   \sum_{n=1}^\infty\frac{x^{2n}}{{2n\choose n}n^3},
   \qquad |x|\le 2,
\label{s3x}
\end{equation}
in terms of periodic zeta functions.  Proposition~\ref{MP} below
shows how  periodic zeta functions may be expressed in terms of
Dirichlet $L$-series. Thus, we can rewrite~(\ref{s3x}) in terms of
$L$-series values, thereby obtaining our result.

\section{Preliminaries}

Let $m$ be a positive integer.  We denote the group of Dirichlet
characters of modulus $m$ by $\hat{U}_m$.  The Dirichlet
$L$-series associated with $\chi\in\hat{U}_m$ is
\[
   L(s,\chi)=\sum_{n=1}^\infty \frac{\chi(n)}{n^s},
   \qquad \R(s)>1.
\]
Similarly, for real $\beta$ we define the {\em periodic zeta
function} (a special case of the Lerch transcendent) by
\[
   \Phi(s,\beta)=\sum_{n=1}^\infty \frac{e^{2\pi i\beta n}}{n^s},
   \qquad \R(s)>1.
\]
Let $\zeta_m=e^{2\pi i/m}$.    Throughout, the sum over a complete
set of residues modulo $m$ is denoted by $\sum_{a\bmod m}$ and the
sum over the positive integer divisors of $m$ is denoted by
$\sum_{d\mid m}$. Thus, Ramanujan's sum is
\[
   c_m(k) = \sum_{\textstyle{\nu\bmod m\atop (\nu,m)=1}} \zeta_m^{\nu k},
\]
and likewise the Gaussian sum attached to $\chi$ is
\[
   \tau(\chi) = \sum_{\nu\bmod m} \chi(\nu)\zeta_m^{\nu}.
\]
Also, $\overline\chi$ denotes the inverse---or equivalently, the
complex conjugate---of the character $\chi$.  Finally, as
customary,  $\mu()$, $\varphi()$, and $\zeta()$ denote the
M\"obius, Euler totient, and Riemann zeta functions, respectively.

Our immediate goal is to represent periodic zeta functions in
terms of $L$-series.  It turns out to be easier to do the reverse
first.  The following result is well known, so we omit the proof.

\begin{lem}\label{L1}
Let $m$ be a positive integer, let $\chi$ be a Dirichlet character
of modulus $m$, and let $L(s,\chi)$ be the associated Dirichlet
$L$-series. Then
\[
   L(s,\chi)=\frac{1}{m}\sum_{a\bmod m}\chi(a)\sum_{b\bmod m}\zeta_m^{-ab}\,\Phi(s,b/m),
   \qquad \R(s)>1.
\]
\end{lem}

\medskip

\begin{lem}\label{L2} Let $a$ and $m$ be positive integers.  Then
\[
   \frac1{\varphi(m)}\sum_{\chi\in\hat{U}_m}\chi(a)\tau(\overline{\chi})L(s,\chi)=
   \frac1{m}\sum_{b\bmod m}\Phi(s,b/m)\,c_m(a-b),
   \qquad \R(s)>1.
\]
\end{lem}
{\em Proof.} First  recall that
\[
   \sum_{\chi\in\hat{U}_m}\overline{\chi}(c)\chi(a)=
   \left\{ \begin{array}{ll}
            \varphi(m) & \;\; \mbox{if $(ac,m)=1$ and $a\equiv c\bmod m$,}\\
                    0  & \;\; \mbox{otherwise.}
           \end{array}\right.
\]
We claim that if $(c,m)=1$, then
\begin{equation}\label{E1}
   \frac{\varphi(m)}{m}\sum_{b\bmod m}\zeta_m^{-bc}\,\Phi(s,b/m)
   =\sum_{\chi\in\hat{U}_m}\chi(c)L(s,\chi);
\end{equation}
for by Lemma \ref{L1},
\begin{eqnarray*}
   \sum_{\chi\in\hat{U}_m}\overline{\chi}(c)\;L(s,\chi)
   &=& \frac1{m}\sum_{a\bmod m}\;\sum_{b\bmod m}\zeta_m^{-ab}\,\Phi(s,b/m)\sum_{\chi\in\hat{U}_m}
       \overline{\chi}(c)\chi(a)\\
   &=& \frac{\varphi(m)}{m}\sum_{b\bmod m}\zeta_m^{-bc}\,\Phi(s,b/m).
\end{eqnarray*}
On the other hand, if $(c,m)>1$, then clearly
\[
   \sum_{\chi\in\hat{U}_m}\overline{\chi}(c)\,L(s,\chi)=0.
\]

We now multiply equation~(\ref{E1}) by $\zeta_m^{ac}$ with
$(a,m)=1$, and then sum over all $c$ modulo $m$, obtaining
\begin{eqnarray*}
   \frac{\varphi(m)}{m}\sum_{b\bmod m}\;\sum_{\textstyle{c\bmod m\atop (c,m)=1}}\zeta_m^{(a-b)c}
   \,\Phi(s,b/m)
   &=& \sum_{c\bmod m}\zeta_m^{ac}\sum_{\chi\in\hat{U}_m}\overline{\chi}(c)L(s,\chi)\\
   &=& \sum_{c\bmod m}\zeta_m^{ac}\sum_{\chi\in\hat{U}_m}\overline{\chi}(c)L(s,\chi)\\
   &=& \sum_{\chi\in\hat{U}_m}\;\sum_{c\bmod m}\overline{\chi}(c)\zeta_m^{ac}\,L(s,\chi)\\
   &=& \sum_{\chi\in\hat{U}_m}\chi(a)\tau(\overline{\chi})L(s,\chi).
\end{eqnarray*}
Rewriting this latter equation in terms of Ramanujan sums
completes the proof. \eop

\medskip

We now state the main proposition of this section.
\begin{prop}\label{MP} Let $a$ and $m$ be coprime positive integers.  Then
\[
   m^s\Phi(s,a/m)=\sum_{d\mid m}\frac{d^s}{\varphi(d)}\sum_{\chi\in\hat{U}_d}
      \chi(a)\tau(\overline{\chi})\,L(s,\chi),\
      \qquad \R(s)>1.
\]
\end{prop}

Before proving Proposition~\ref{MP}, we state and prove two
lemmata which are used in the proof of Proposition~\ref{MP}.

\begin{lem}\label{LM} Let $f:\Z\to\C$ be multiplicative and such
that for all positive integers $m$,
\[
   F(m):=\sum_{d|m}\mu^2(d)f(d)
\]
is non-zero. Furthermore, let
\[
   g(m):=\sum_{d|m}\frac{\mu(d)}{F(d)}.
\]
Then for all positive integers $k$ and $m$ such that $k$ divides
$m$,
\[
   \sum_{\textstyle{d|m \atop k|d}}\mu^2(d)f(d)=F(m)\mu^2(k)g(k).
\]
In particular,
\[
   \sum_{\textstyle{d|m\atop k|d}}\frac{\mu^2(d)}{\varphi(d)}
   =\frac{m}{\varphi(m)}\frac{\mu^2(k)}{k}.
\]
\end{lem}

\noindent{\it Proof.} First, let us define $F(x)=f(x)=0$ if $x$ is
not an integer.  Next, observe that $F(p^a)=F(p)$ for all positive
primes $p$ and positive integers $a$.  We may write $m$ as
$\prod_pp^{a_p}$ and $k$ as $\prod_pp^{b_p}$ where $p$ ranges over
all positive primes and $a_p$ and $b_p$ are non-negative integers
with $b_p\leq a_p$.  Since $F$ is multiplicative, we have
\begin{eqnarray}
   \sum_{\textstyle{d|m \atop k|d}}\mu^2(d)f(d)
   &=&
   \prod_{p|m}\bigg(\sum_{\nu_p=b_p}^{a_p}\mu^2(p^{\nu_p})f(p^{\nu_p})\bigg)\nonumber\\
   &=& \prod_{p|m}\left(F(p^{a_p})-F(p^{b_p-1})\right)\nonumber\\
   &=& \prod_{p|m}F(p^{a_p})\left(1-\frac{F(p^{b_p-1})}{F(p^{a_p})}\right)\nonumber\\
   &=& F(m)\prod_{p|m}\left(1-\frac{F(p^{b_p-1})}{F(p^{a_p})}\right).
\label{Fprod}
\end{eqnarray}
Notice that the final product in~(\ref{Fprod}) vanishes if any
$b_p\geq 2$, for then $F(p^{b_p-1})= F(p)= F(p^{a_p})$.  Hence if
$k$ is not square-free, then the lemma is trivially true as both
sides are equal to $0$.  Therefore, we may assume henceforth that
$k$ is square-free.  Now if $b_p=0$, then
$1-F(p^{b_p-1})/F(p^{a_p})=1$, and thus (under the assumption that
$k$ is square-free), we may restrict the final product
in~(\ref{Fprod}) to primes $p$ for which $b_p=1$.  This yields
\begin{eqnarray*}
   \sum_{\textstyle{d|m \atop k|d}}\mu^2(d)f(d)
   &=& F(m)\prod_{p|k}\left(1-\frac{1}{F(p^{a_p})}\right)
   =F(m)\prod_{p|k}\left(1-\frac{1}{F(p)}\right)\\
   &=& F(m)\sum_{d|k}\frac{\mu(d)}{F(d)}\\
   &=& F(m)g(k).
\end{eqnarray*}
Thus, in general, we have
\[
   \sum_{\textstyle{d|m \atop k|d}}\mu^2(d)f(d)=F(m)\mu^2(k)g(k).
\]
The special case is obtained by taking $F(m)=m/\varphi(m)$, so
that if $k$ is square-free, then
\[
   g(k) = \sum_{d|k}\frac{\mu(d)\varphi(d)}{d} =
   \prod_{p|k}\left(1-\frac{\varphi(p)}{p}\right)
   =\prod_{p|k}\frac1p=\frac{1}{k}.
\]
This completes the proof of Lemma~\ref{L1}. \eop

\begin{lem}\label{LB} Let $m$ be a positive integer and let $\beta$ be any real number.
Then
\[
   \sum_{\textstyle{n=1\atop (n,m)=1}}^\infty \frac{e^{2\pi i\beta n}}{n^s}
       =\sum_{d\mid m}\frac{\mu(d)}{d^s}\,\Phi(s,\beta d).
\]
\end{lem}
{\em Proof.} Let $x$ be any complex number with $|x|\leq 1$. Then
\begin{eqnarray*}
   \sum_{\textstyle{n=1\atop (n,m)=1}}^\infty\frac{x^n}{n^s}
   &=&\sum_{n=1}^\infty\frac{x^n}{n^s}\sum_{d\mid (n,m)}\mu(d)
   =\sum_{n=1}^\infty\frac{x^n}{n^s}\sum_{\textstyle{d\mid n\atop d\mid m}}\mu(d)
   =\sum_{d\mid m}\mu(d)\sum_{k=1}^\infty\frac{x^{kd}}{(kd)^s}\\
   &=&\sum_{d\mid m}\frac{\mu(d)}{d^s}\sum_{k=1}^\infty\frac{x^{kd}}{k^s}.
\end{eqnarray*}
Replacing $x$ by $e^{2\pi i\beta}$ completes the proof. \eop

\medskip

\noindent {\em Proof of Proposition~\ref{MP}.} First, recall (see
eg.~\cite[p.\ 238]{HW}) that Ramanujan's sum has the explicit
representation
\[
   c_m(k)=\varphi(m)\frac{\mu(m/(m,k))}{\varphi(m/(m,k))}.
\]
Hence, we have
\begin{eqnarray}
   &&\frac{1}{m}\sum_{b\bmod m}\Phi(s,b/m)\,c_m(a-b)\nonumber\\
   &=& \frac{1}{m}\sum_{b\bmod m}\Phi(s,b/m)
       \varphi(m)\frac{\mu(m/(m,a-b))}{\varphi(m/(m,a-b))}\nonumber\\
   &=&\frac{\varphi(m)}{m}\sum_{d\mid m}\sum_{\textstyle{b\bmod m\atop
       (a-b,m)=m/d}} \Phi(s,b/m)\frac{\mu(d)}{\varphi(d)}\nonumber\\
   &=&\frac{\varphi(m)}{m}\sum_{d\mid m}\frac{\mu(d)}{\varphi(d)}
      \sum_{\textstyle{\nu\bmod d\atop (\nu,d)=1}}\Phi\bigg(s,\frac{a+m\nu/d}{m}\bigg)\nonumber\\
   &=&\frac{\varphi(m)}{m}\sum_{d\mid m}\frac{\mu(d)}{\varphi(d)}\sum_{n=1}^\infty
      n^{-s}\sum_{\textstyle{\nu\bmod d\atop (\nu,d)=1}}\zeta_m^{(a+m\nu/d)n}\nonumber\\
   &=&\frac{\varphi(m)}{m}\sum_{d\mid m}\frac{\mu(d)}{\varphi(d)}\sum_{n=1}^\infty
      n^{-s}\zeta_m^{an}c_d(n)\nonumber\\
   &=&\frac{\varphi(m)}{m}\sum_{d\mid
        m}\frac{\mu(d)}{\varphi(d)}\sum_{n=1}^\infty
      n^{-s}\zeta_m^{an}\varphi(d)\frac{\mu(d/(n,d))}{\varphi(d/(n,d))}\nonumber\\
   &=&\frac{\varphi(m)}{m}\sum_{d\mid m}\mu(d)\sum_{f\mid
      d}\frac{\mu(f)}{\varphi(f)} \sum_{\textstyle{n=1\atop
      (n,f)=1}}^\infty(nd/f)^{-s}\zeta_{fm/d}^{an}\nonumber\\
   &=&\frac{\varphi(m)}{m}\sum_{d\mid m}\mu(d)\sum_{f\mid
      d}\frac{\mu(f)}{\varphi(f)} \bigg(\frac{d}{f}\bigg)^{-s}\sum_{\textstyle{n=1\atop
     (n,f)=1}}^\infty n^{-s}\zeta_{fm/d}^{an}.
\label{saveit}
\end{eqnarray}
But by Lemma~\ref{LB} the final expression in~(\ref{saveit}) can
be rewritten as
\begin{equation}
   \frac{\varphi(m)}{m}\sum_{d\mid m}\mu(d)\sum_{f\mid
    d}\frac{\mu(f)}{\varphi(f)} \left(\frac{d}{f}\right)^{-s}\sum_{\delta\mid
    f}\delta^{-s}\mu(\delta)\Phi\bigg(s,\frac{ad\delta}{fm}\bigg).
\label{transformit}
\end{equation}
Now transform~(\ref{transformit}) by changing the variable $f$ to
$d/f$, then letting $k=f\delta$ (noticing that the only non-zero
terms occur when $d$ is square-free), then observing that
$\sum_{f\mid k}\varphi(f)=k$, and finally replacing $d$ by $kd$.
Thus, from~(\ref{saveit}) and~(\ref{transformit}),
\begin{eqnarray}
   &&\frac{1}{m}\sum_{b\bmod m}\Phi(s,b/m)\,c_m(a-b)\nonumber\\
   &=&\frac{\varphi(m)}{m}\sum_{d\mid m}\mu(d)\sum_{f\mid
      d}\frac{\mu(d/f)}{\varphi(d/f)} f^{-s}\sum_{\delta\mid
      (d/f)}\delta^{-s}\mu(\delta)\Phi(s,af\delta/m)\nonumber\\
   &=&\frac{\varphi(m)}{m}\sum_{d\mid m}\sum_{k\mid
      d}\frac{\mu^2(d)}{\varphi(d)}
      k^{-s}\mu(k)\Phi(s,ak/m)\sum_{f\mid k}\varphi(f)\nonumber\\
   &=&\frac{\varphi(m)}{m}\sum_{k\mid m}k^{1-s}\mu(k)\Phi(s,ak/m)
    \sum_{\textstyle{d\mid m\atop k\mid
    d}}\frac{\mu^2(d)}{\varphi(d)}.
\label{DoSomeMore}
\end{eqnarray}
By applying Lemma~\ref{LM} to~(\ref{DoSomeMore}) and then
replacing $m/k$ by $d$, we find that
\begin{eqnarray}
   \frac{1}{m}\sum_{b\bmod m}\Phi(s,b/m)\,c_m(a-b)
   &=&\sum_{k\mid m}k^{-s}\mu(k)\Phi(s,ak/m)\nonumber\\
   &=&\frac{1}{m^s}\sum_{d\mid m}d^s\mu(m/d)\Phi(s,a/d).
\label{AlmostThere}
\end{eqnarray}
Hence by~(\ref{AlmostThere}) and Lemma~\ref{L2}, we see that
\begin{eqnarray*}
   \frac{1}{m^s}\sum_{d\mid m}d^s\mu(m/d)\Phi(s,a/d)
   &=& \frac{1}{m}\sum_{b\bmod m}\Phi(s,b/m)\,c_m(a-b)\\
   &=&\frac{1}{\varphi(m)}\sum_{\chi\in\hat{U}_m}\chi(a)\tau(\overline{\chi})L(s,\chi).
\end{eqnarray*}
An application of M\"obius inversion now completes the proof. \eop

\section{Main Results}
We are now in a position to derive our class number formula. To
this end, for $|x|\le 2$ and $2\le k\in\Z$, put
\[
   s(k,x):=\sum_{n=1}^\infty \frac{x^{2n}}{{2n\choose n}n^k}.
\]
Let $0<\theta<\pi$ and $x=2\sin\theta/2$.  Then~\cite[p.\ 61
(2)]{GR} $2s(2,x)=\theta^2$ and by formula (2.7) of~\cite{Z},
\begin{eqnarray}
   \theta^2\log(2\sin \theta/2)
   &=& 2\zeta(3)+
      \sum_{n=1}^\infty \frac{(2\sin \theta/2)^{2n}}{{2n\choose n}n^3}
      \nonumber\\
   &-& 2\theta\,\I \;\Phi(2,\theta/2\pi)-2\R\;\Phi(3,\theta/2\pi),
\label{zuck}
\end{eqnarray}
where $\R$ and $\I$ denote the real and imaginary parts of a
complex number, respectively. Now substitute $\theta=2\pi a/m$
with $(a,m)=1$ and $0<a<m/2$ in~(\ref{zuck}) to obtain
\begin{eqnarray}\label{E2}
   \log\left(2\sin \frac{\pi a}{m}\right)
   &=&\frac{m^2}{2\pi^2a^2}\zeta(3)+\frac{m^2}{4\pi^2a^2}
   \sum_{n=1}^\infty \frac{(2\sin \pi a/m)^{2n}}{{2n\choose
   n}n^3}\nonumber\\
   &-&\frac{m}{\pi a}\I \; \Phi\left(2,\frac{a}{m}\right)
    -\frac{m^2}{2\pi^2a^2}\R \;\Phi\left(3,\frac{a}{m}\right).
\end{eqnarray}

In our main result, character sums of consecutive integer powers
arise, and it is convenient to fix some notation for these.

\begin{Def}\label{Bchi} Let $m$ be a positive integer.  If $\chi$ is a
Dirichlet character of modulus $m$ and $j$ is any integer, put
\begin{equation}
   \mathcal{B}_j(\chi) := \sum_{0<a< m/2} a^j\chi(a).
\label{mathcalB}
\end{equation}
\end{Def}
We now state and prove our main result.

\begin{thm}\label{T} Let $m$ be a positive integer, let
$\psi$ be an even primitive character of modulus $m$, and let
$\mathcal{B}_j$ be as in~(\ref{mathcalB}).  Then
\begin{eqnarray*}
   L(1,\psi)
   &=&\frac{2\tau(\psi)}{\pi i m^2}\sum_{d\mid m}\frac{d^2}{\varphi(d)}
    \sum_{\textstyle{\chi\in\hat{U}_d\atop \chi\,\mbox{\footnotesize odd}}}
    \mathcal{B}_{-1}(\chi\overline{\psi})\tau(\overline{\chi})L(2,\chi)\\
   &+&\frac{\tau(\psi)}{\pi^2 m^2}\sum_{d\mid m}\frac{d^3}{\varphi(d)}
    \sum_{\textstyle{\chi\in\hat{U}_d\atop \chi\,\mbox{\footnotesize even}}}
    \mathcal{B}_{-2}(\chi\overline{\psi})\tau(\overline{\chi})\,L(3,\chi)
  -\frac{m\tau(\psi)}{\pi^2}\mathcal{B}_{-2}(\overline{\psi})\zeta(3)\\
  &-&\frac{m\tau(\psi)}{2\pi^2}\sum_{n=1}^\infty
  \frac{1}{{2n\choose n}n^3}\sum_{0<a<m/2}\frac{\overline{\psi}(a)}{a^2}
  \left(2\sin\frac{\pi a}{m}\right)^{2n}.
\end{eqnarray*}
\end{thm}
{\em Proof.} We start with~(\ref{E2}) and write $\I \;
\Phi(2,a/m)$ and $\R \; \Phi(3,a/m)$ in terms of $L$-series via
Proposition \ref{MP}; but first observe that
\begin{eqnarray*}
   &&\I \bigg(\sum_{\chi\bmod
   d}\chi(a)\tau(\overline{\chi})L(2,\chi)\bigg)\\
   &=&\frac1{2i}\sum_{\chi\in\hat{U}_d}\left(\chi(a)\tau(\overline{\chi})L(2,\chi)
     -\overline{\chi(a)\tau(\overline{\chi})\,L(2,\chi)}\right)\\
   &=&\frac1{2i}\sum_{\chi\in\hat{U}_d}\left(\chi(a)\tau(\overline{\chi})L(2,\chi)
       -\chi(-1)\overline{\chi}(a)\tau(\chi)L(2,\overline{\chi})\right),
\end{eqnarray*}
since $\overline{\tau(\chi)}=\chi(-1)\tau(\overline{\chi})$. Now
split the sum over the two terms and in the second sum replace
$\chi$ by $\overline{\chi}$.  The even characters cancel and we
obtain
\[
   \I \bigg(\sum_{\chi\in\hat{U}_d}
      \chi(a)\tau(\overline{\chi})\,L(2,\chi)\bigg)
   = \frac{1}{i}\sum_{\textstyle{\chi\in\hat{U}_d\atop \chi(-1)=-1}}
      \chi(a)\tau(\overline{\chi})\,L(2,\chi).
\]
Similarly, we see that
\[
   \R \bigg(\sum_{\chi\in\hat{U}_d}
      \chi(a)\tau(\overline{\chi})\,L(3,\chi)\bigg)
   = \sum_{\textstyle{\chi\in\hat{U}_d\atop \chi(-1)=1}}
      \chi(a)\tau(\overline{\chi})\,L(3,\chi).
\]
Thus by~(\ref{E2}) and Proposition~\ref{MP},
\begin{eqnarray}
   \log\left(2\sin\frac{\pi a}{m}\right)
   &=&\frac{m^2}{2\pi^2a^2}\zeta(3)
     +\frac{m^2}{4\pi^2a^2}
     \sum_{n=1}^\infty\frac{\left(2\sin\pi a/m\right)^{2n}}{{2n\choose
     n}n^3}\nonumber\\
   &-& \frac{1}{\pi i ma}
      \sum_{d\mid m}\frac{d^2}{\varphi(d)}
      \sum_{\textstyle{\chi\in\hat{U}_d\atop \chi\,\mbox{\footnotesize
      odd}}}\chi(a)\tau(\overline{\chi})\,L(2,\chi)\nonumber\\
   &-&\frac{1}{2\pi^2ma^2}\sum_{d\mid m}\frac{d^3}{\varphi(d)}
      \sum_{\textstyle{\chi\in\hat{U}_d\atop \chi\,\mbox{\footnotesize
      even}}}\chi(a)\tau(\overline{\chi})\,L(3,\chi).
\label{logformula}
\end{eqnarray}
Next, recall (see eg.~\cite[p.\ 336]{BS}) that if $\psi$ is an
even primitive character of modulus $m$, then
\begin{eqnarray}
   L(1,\psi)
   &=&-\frac{\tau(\psi)}{m}\sum_{a=1}^{m-1}
       \overline{\psi}(a)\log\left(2\sin\frac{\pi a}{m}\right)\nonumber\\
   &=&-\frac{2\tau(\psi)}{m}\sum_{0<a<m/2}
   \overline{\psi}(a)\log\left(2\sin\frac{\pi a}{m}\right).
\label{gaussformula}
\end{eqnarray}
Substituting~(\ref{gaussformula}) into~(\ref{logformula})
completes the proof. \eop

\medskip

Let $D$ be a (positive fundamental) discriminant of a real
quadratic number field. Let $h(D)$ denote its class number,
$\varepsilon=\varepsilon_D$ its fundamental unit $>1$, and
$\chi_D=(D/\cdot)$, the Kronecker symbol, i.e. the Dirichlet
character associated with the quadratic field of discriminant $D$.
Then by Dirichlet (see eg.~\cite[p.\ 343]{BS}), we know that
$$2h(D)\log\eps_D=\sqrt{D}\;L(1,\chi_D).$$  Hence by
Theorem~\ref{T}, using the fact that $\tau(\chi_D)=\sqrt{D}$, we
obtain the following class number formula.

\begin{cor}[Class Number Formula] \label{C}
\begin{eqnarray*}
   h(D)\log\eps_D
   &=&\frac{1}{\pi Di}\sum_{d\mid D}\frac{d^2}{\varphi(d)}
   \sum_{\textstyle{\chi\in \hat{U}_d\atop \chi\,\mbox{\footnotesize odd}}}\mathcal{B}_{-1}(\chi\chi_D)
   \tau(\overline{\chi})\,L(2,\chi)\\
   &+&\frac{1}{2\pi^2 D}\sum_{d\mid D}\frac{d^3}{\varphi(d)}
   \sum_{\textstyle{\chi\in \hat{U}_d\atop \chi\,\mbox{\footnotesize even}}}\mathcal{B}_{-2}(\chi\chi_D)
   \tau(\overline{\chi})\,L(3,\chi)\\
   &-&\frac{D^2}{2\pi^2}\mathcal{B}_{-2}(\chi_D)\zeta(3)\\
   &-&\frac{D^2}{4\pi^2}\sum_{n=1}^\infty
   \frac{1}{{2n\choose n}n^3}\sum_{0<a<D/2}\frac{\chi_D(a)}{a^2}\left(2\sin\frac{\pi a}{D}\right)^{2n}.
\end{eqnarray*}
\end{cor}

\subsection{A Computation}

As an amusing conclusion, we now show how to use our class number
formula to compute $h(5)$, the class number of the quadratic field
$\Q(\sqrt{5}).$  Since the discriminant $D=5$, the only relevant
moduli of characters are $m=1$ and $m=5$. For $m=1$, the unique
character is the even constant character $1$.  For $m=5$, we have
four characters determined by the homomorphisms from
$(\Z/5\Z)^\times$ into $\C^\times$, namely $\chi_\nu$ for
$\nu=0,1,2,3$ determined by $\chi_\nu(2)=i^\nu$. Notice that
$\overline{\chi_1}= \chi_3$ and that $\chi_2=(5/\cdot)=\chi_5$,
the Kronecker character modulo $5$.

By Corollary~\ref{C}, we have $h(5)=(A+B+C+S)/\log \eps_5$, where
\begin{eqnarray*}
   A &=&\frac{5}{4\pi i}\left(\mathcal{B}_{-1}(\chi_3)\tau(\chi_3)L(2,\chi_1)+
        \mathcal{B}_{-1}(\chi_1)\tau(\chi_1)L(2,\chi_3)\right),\\
   B &=& \frac1{10\pi^2}\bigg(\mathcal{B}_{-2}(\chi_5)\tau(1)\zeta(3)\\
     && \qquad + \frac{125}{4}\big(\mathcal{B}_{-2}(\chi_5)\tau(\chi_0)L(3,\chi_0)
        +\mathcal{B}_{-2}(\chi_0)\tau(\chi_5)L(3,\chi_5)\big)\bigg)\\
   C &=& -\frac{25}{2\pi^2}\mathcal{B}_{-2}(\chi_5)\zeta(3)\\
   S &=& -\frac{25}{4\pi^2}\sum_{n=1}^\infty\frac1{{2n\choose n}n^3}
         \left(\chi_5(1)(2\sin\pi/5)^{2n}+\tfrac14\chi_5(2)(2\sin 2\pi/5)^{2n}\right),
\end{eqnarray*}
and
\begin{eqnarray*}
   \mathcal{B}_{-1}(\chi_1) &=& \chi_1(1)+\tfrac12\chi_1(2)=1+\tfrac{1}{2}i\\
   \mathcal{B}_{-1}(\chi_3) & =& 1-\tfrac1{2}i\\
   \mathcal{B}_{-2}(\chi_0) &=& 1+\tfrac1{4}=\tfrac5{4}\\
   \mathcal{B}_{-2}(\chi_5) &=& 1-\tfrac1{4}=\tfrac3{4}\\
   \tau(1) &=& 1\\
   \tau(\chi_0) &=& -1\\
   \tau(\chi_5) &=& \sqrt{5}\\
   \tau(\chi_1) &=& \zeta_5+i\zeta_5^2-i\zeta_5^3-\zeta_5^4=
         \left(i+\frac{1-\sqrt{5}}{2}\right)\sqrt{\frac{5+\sqrt{5}}{2}}\\
   \tau(\chi_3) &=& \zeta_5-i\zeta_5^2+i\zeta_5^3-\zeta_5^4=
         \left(i+\frac{\sqrt{5}-1}{2}\right)\sqrt{\frac{5+\sqrt{5}}{2}}\\
   L(3,\chi_0) &=& \left(1-5^{-3}\right)\zeta(3)=\frac{124}{125}\zeta(3).
\end{eqnarray*}

In order to evaluate $L(s,\chi_\nu)$ for $\nu=0,1,2,3$ and $s=2,3$
we write
\[
   L(s,\chi_\nu)=5^{-s}\sum_{r=1}^4\chi_\nu(r)\zeta(s,r/5),
\]
where $\zeta(s,a)=\sum_{n=0}^\infty (n+a)^{-s}$ is the Hurwitz
zeta function.  Hence to evaluate these $L$-series, it suffices to
evaluate the Hurwitz zeta functions. The following table gives the
appropriate approximations.

\begin{center}
\begin{tabular}{|l||l|l|}\hline
   $r$ & $\zeta(2,r/5)$ & $\zeta(3,r/5)$\\ \hline\hline
   $1$ &  26.26737720... & 125.73901805...\\
   \hline
   $2$ & 7.27535659...& 16.1195643...\\
   \hline
   $3$ & 3.63620967...& 4.98141576... \\
   \hline
   $4$ & 2.29947413...& 2.21505785...\\
   \hline
\end{tabular}
\end{center}
Thus, we find that
\begin{eqnarray*}
   L(2,\chi_1) &=& 0.95871612... + (0.14556587...)i\\
   L(2,\chi_3) &=& 0.95871612... - (0.14556587...)i\\
   L(3,\chi_5) &=& 0.85482476...
\end{eqnarray*}
Furthermore, $\zeta(3)=1.20205690...$, so that
\begin{eqnarray*}
   A &=&  1.24907310...\\
   B &=&  0.48248793...\\
   C &=& -1.14181713...\\
   S &=& -0.10853146...
\end{eqnarray*}
Finally, the fundamental unit $$\eps_5=\frac{1+\sqrt{5}}{2},$$
(which generally can be computed efficiently by continued
fractions).

From all of this we obtain $h(5)=1.0000000...+(0.0000000...)i$, whence $h(5)=1$.


{\bf Address of authors:}\\
\indent\indent        Department of Mathematics and Statistics\\
\indent\indent        University of Maine \\
\indent\indent        Orono, Maine 04469-5752\\
\indent\indent        U.S.A.
\\
\\
{\bf e-Mail Addresses:}\\ \indent\indent bradley@math.umaine.edu\\
\indent\indent ozluk@math.umaine.edu\\  \indent\indent
snyder@math.umaine.edu

\begin{thebibliography}{99}

\bibitem{B}
Bradley, D., A class of series acceleration formulae for Catalan's
constant, {\it The Ramanujan J.}, {\bf 3} (1999), no.~2, 159-173.
[MR 1703281] (2000f:11163) {\tt http://arxiv.org/abs/0706.0356} 

\bibitem{BS}
Borevich, Z. and Shafarevich, I., {\it Number Theory},
Academic Press, New York and London, (1966).

\bibitem{GR}
Gradshteyn, I.\ S.\ and Ryzhik, I.\ M., {\it Table of Integrals,
Series, and Products} (5th ed.), Academic Press, Boston, 1994.

\bibitem{HW}
Hardy, G.\ H.\ and Wright, E.\ M., {\it An Introduction to the
Theory of Numbers} (5th ed.), Clarendon Press, Oxford, 1979.


\bibitem{Z}
Zucker, I.\ J., On the series $\sum_{k=1}^\infty {2k\choose
k}^{-1}k^{-n}$ and related sums, {\it J.\ Number Theory}, {\bf 20}
(1985), 92-102.

\end{thebibliography}
\end{document}